\documentclass[draft]{amsart}

\usepackage{amssymb}
\usepackage{url}
\usepackage{amscd}

\theoremstyle{plain}
\newtheorem{theorem}{Theorem}

\theoremstyle{remark}

\begin{document}

\date{} 

\title{Minimal hypersurfaces with zero Gauss-Kronecker curvature}

\author{T. Hasanis,\ A. Savas-Halilaj and\ T.Vlachos}

\address{Department of Mathematics\\
University of Ioannina\\
45110-Ioannina\\
Greece}

\email{thasanis@cc.uoi.gr, me00499@cc.uoi.gr, tvlachos@cc.uoi.gr}


\subjclass[2000]{53C42}

\begin{abstract}
We investigate complete minimal hypersurfaces in the Euclidean space $%
\ {R}^{4}$, with Gauss-Kronecker curvature identically zero. We prove
that, if $f:M^{3}\rightarrow \ {R}^{4}$ is a complete minimal
hypersurface with Gauss-Kronecker curvature identically zero, nowhere
vanishing second fundamental form and scalar curvature bounded from below,
then $f(M^{3})$ splits as a Euclidean product $L^{2}\times 
\ {R}$, where $L^{2}$ is a complete minimal surface in $\ {R}^{3}$
with Gaussian curvature bounded from below.
\end{abstract}

\maketitle


\newcommand\sfrac[2]{{#1/#2}}

\newcommand\cont{\operatorname{cont}}
\newcommand\diff{\operatorname{diff}}


\section{Introduction}
A classical result of Beez-Killing states that a hypersurface $M^{n}$ in the
Euclidean space $R^{n+1}$ is rigid if the rank of the Gauss map is at least
3. Dajczer and Gromoll in \cite{DG} developed a powerful method, the so
called ``\textit{Gauss parametrization}'', which has interesting
applications in the study of rigidity problems in the case where the rank is
at least 2. The essential point of this method is that it provides a
parametrization for every hypersurface with second fundamental form of
constant nullity by inverting the Gauss map. The local rigidity of minimal
hypersurfaces with nullity $n-2$ is well understood. In fact, Dajczer and
Gromoll \cite{DG}, show that such hypersurfaces allow locally an
one-parameter family of isometric deformations, the so called associated
family. Hence, the rigidity of minimal hypersurfaces make sense under global
assumptions. Among others, Dajczer and Gromoll in \cite{DG} proved the
following global rigidity result: \textit{If }$M^{n}$\textit{\ }$\left(
n\geq 4\right) $\textit{\ is a complete Riemannian manifold which does not
have }$R^{n-3}$\textit{\ as a factor, then every minimal immersion }$%
f:M^{n}\rightarrow R^{n+1}$\textit{\ is rigid as a
minimal submanifold of }$R^{n+p}$\textit{\ for any }$%
p\geq 1$\textit{. }This is proved by establishing an interesting criterion
for a complete hypersurface $M^{n}$ in $R^{n+1}$ to split as a Euclidean
product $M^{n}=L^{3}\times R^{n-3}$. However, the case of dimension $n=3$ is
left open.

The aim of the present paper is to fill in the gap for $n=3$, by studying
complete minimal hypersurfaces in $R^{4}$ with Gauss-Kronecker curvature
identically zero. One can easily construct complete minimal hypersurfaces
with Gauss-Kronecker curvature identically zero in $R^{4}$ just by erecting
cylinders over complete minimal surfaces in $R^{3}$. In fact, let $%
h:M^{2}\rightarrow R^{3}\hookrightarrow R^{4}$ be\ a complete\ minimal
surface. Denote by $\lambda $, $-\lambda $ its principal curvatures and by $%
a $ a unit vector in $\mathbb{R}^{4}$ normal to $R^{3}$. Then the cylinder $%
f:M^{3}=M^{2}\times R\rightarrow R^{4}$, $f\left( p,t\right) =h\left(
p\right) +ta$, is a complete minimal hypersurface in $R^{4}$ with principal
curvatures $k_{1}=\lambda $, $k_{2}=0$, $k_{3}=-\lambda $. The scalar
curvature of $M^{3}$ is equal to the Gaussian curvature of $M^{2}$.\ Observe
that, the cylinder over the helicoid in $R^{3}$ gives a non-totally geodesic
complete minimal hypersurface in $R^{4}$, with Gauss-Kronecker curvature
identically zero, three distinct principal curvatures and scalar curvature
bounded from below.\textit{\ }Also, we note that there are complete minimal
surfaces in $R^{3}$ having their Gaussian curvature not bounded from below 
\cite{Na}. So, erecting the cylinder over those surfaces, one can produce
complete minimal hypersurfaces in $R^{4}$ with Gauss-Kronecker curvature
identically zero and scalar curvature not bounded from below. However, to
the best of our knowledge, \ the cylinders are the only known examples of
complete minimal hypersurfaces with Gauss-Kronecker identically zero in $%
R^{4}$. This led us to the following\textit{\medskip }

\noindent \textbf{Question: }\textit{It is true that any complete minimal
hypersurface with vanishing Gauss-Kronecker curvature in }$R^{4}$\textit{\
is a cylinder over a minimal surface in }$R^{3}$\textit{?\medskip }

We shall give here a partial answer to this question under some additional
assumptions on the scalar curvature. In particular, we prove the following

\begin{theorem}
Let $M^{3}$ be an oriented, 3-dimensional, complete
Riemannian manifold and $f:M^{3}\rightarrow R^{4}$ a minimal
isometric immersion with Gauss-Kronecker curvature identically zero and
nowhere vanishing second fundamental form. If the scalar curvature is
bounded from below, then $f\left( M^{3}\right)$ splits as a Euclidean
product$\ L^{2}\times R$, where $L^{2}$ is a complete
minimal surface in $R^{3}$ with Gaussian curvature bounded from
below.
\end{theorem}

\noindent \textbf{Remark}:\textbf{\ }It is clear that the cylinders in
Theorem 1 are not rigid.\medskip

Recently Cheng \cite{C} proved that complete minimal hypersurfaces in $R^{4}$
with scalar curvature bounded from below and constant Gauss-Kronecker
curvature, have identically zero Gauss-Kronecker curvature. Also he gave an
example of a minimal hypersurface in $R^{4}$ with Gauss-Kronecker curvature
identically zero which is in fact a cone shaped hypersurface over the
Clifford torus in $S^{3}$ and certainly is not complete. Motivated by these
facts, we consider complete minimal hypersurfaces in $R^{4}$ with constant
Gauss-Kronecker curvature and prove that this constant must be zero without
any assumption on the scalar curvature. The key of the proof is the powerful
Principal Curvature Theorem proved by Smyth and Xavier in \cite{SX}. More
precisely, we prove the following

\begin{theorem}
Let $M^{3}$ be an oriented, 3-dimensional, complete Riemannian
manifold and $f:M^{3}\rightarrow R^{4}$ a minimal isometric
immersion with constant Gauss-Kronecker curvature. Then the Gauss-Kronecker
curvature is identically zero.
\end{theorem}


\section{Preliminaries}
Let $f:M^{3}\rightarrow R^{4}$ be an oriented minimal hypersurface equipped
with the induced metric and unit normal vector field $\xi $ along $f$.
Denote by $A$ the shape operator associated with $\xi $ and by $k_{1}\geq
k_{2}\geq k_{3}$ the principal curvatures. The Gauss-Kronecker curvature $K$
and the scalar curvature $\tau $\ are given by 
\begin{equation*}
K=\det A=k_{1}k_{2}k_{3}\text{, }\tau =k_{1}k_{2}+k_{1}k_{3}+k_{2}k_{3}.
\end{equation*}
Assume now that the second fundamental form is nowhere vanishing. Then the
principal curvatures satisfies the relation $k_{1}=\lambda
>k_{2}=0>k_{3}=-\lambda $, where $\lambda $ is a smooth positive function on 
$M^{3}$. We can choose locally an orthonormal frame field $\left\{
e_{1},e_{2},e_{3}\right\} $ of principal directions corresponding to $%
\lambda ,0,-\lambda $. Let $\left\{ \omega _{1},\omega _{2},\omega
_{3}\right\} $ and $\left\{ \omega _{ij}\right\} $, $i,j\in \{1,2,3\}$, be
the corresponding dual and the\ connection forms. Throughout this paper we
make the following convection for indices 
\begin{equation*}
1\leq i,j,k,\ldots \leq 3,
\end{equation*}
and adopt the method of moving frames. The structure equations are 
\begin{eqnarray*}
d\omega _{i} &=&\sum_{j}\omega _{ij}\wedge \omega _{j},\quad \omega
_{ij}+\omega _{ji}=0, \\
d\omega _{ij} &=&\sum_{l}\omega _{il}\wedge \omega _{lj}-k_{i}k_{j}\omega
_{i}\wedge \omega _{j},\ i\neq j.
\end{eqnarray*}
Consider the functions 
\begin{equation*}
u:=\omega _{12}\left( e_{3}\right) ,\quad v:=e_{2}\left( \log \lambda
\right) ,
\end{equation*}
which will play an crucial role in the proof of Theorem 1. From the
structural equations, and the Codazzi equations, 
\begin{equation*}
\begin{array}{l}
e_{i}\left( k_{j}\right) =\left( k_{i}-k_{j}\right) \omega _{ij}\left(
e_{j}\right) ,\ i\neq j,\medskip  \\ 
\left( k_{1}-k_{2}\right) \omega _{12}\left( e_{3}\right) =\left(
k_{2}-k_{3}\right) \omega _{23}\left( e_{1}\right) =\left(
k_{1}-k_{3}\right) \omega _{13}\left( e_{2}\right) ,
\end{array}
\end{equation*}
we easily get 
\begin{equation*}
\begin{array}{lll}
\omega _{12}\left( e_{1}\right) =v, & \omega _{13}\left( e_{1}\right) =\frac{%
1}{2}e_{3}\left( \log \lambda \right) , & \omega _{23}\left( e_{1}\right)
=u,\medskip  \\ 
\omega _{12}\left( e_{2}\right) =0, & \omega _{13}\left( e_{2}\right) =\frac{%
1}{2}u, & \omega _{23}\left( e_{2}\right) =0,\medskip  \\ 
\omega _{12}\left( e_{3}\right) =u, & \omega _{13}\left( e_{3}\right) =-%
\frac{1}{2}e_{1}\left( \log \lambda \right) , & \omega _{23}\left(
e_{3}\right) =-v,
\end{array}
\end{equation*}
and 
\begin{equation}
\left\{ 
\begin{array}{l}
e_{2}\left( v\right) =v^{2}-u^{2},\medskip  \\ 
e_{1}\left( u\right) =e_{3}\left( v\right) ,\medskip  \\ 
e_{2}\left( u\right) =2uv,\medskip  \\ 
e_{3}\left( u\right) =-e_{1}\left( v\right) ,\medskip  \\ 
e_{1}e_{1}\left( \log \lambda \right) +e_{3}e_{3}\left( \log \lambda \right)
-\frac{1}{2}\left( e_{1}\left( \log \lambda \right) \right) ^{2}\smallskip 
\\ 
\ -\frac{1}{2}\left( e_{3}\left( \log \lambda \right) \right)
^{2}-2v^{2}-4u^{2}+2\lambda ^{2}=0.
\end{array}
\right. 
\end{equation}
Furthermore, the\ above equations yield 
\begin{equation*}
\lbrack e_{1},e_{2}]=-ve_{1}+\frac{1}{2}ue_{3},
\end{equation*}
\begin{equation}
\lbrack e_{1},e_{3}]=-\frac{1}{2}e_{3}\left( \log \lambda \right)
e_{1}-2ue_{2}+\frac{1}{2}e_{1}\left( \log \lambda \right) e_{3},
\end{equation}
\begin{equation*}
\lbrack e_{2},e_{3}]=\frac{1}{2}ue_{1}+ve_{3}.
\end{equation*}


\section{Proofs}
We shall use in the proof of Theorem 1 a result due to S.Y Cheng and S.T Yau 
\cite{CY} that we recall in the following lemma. For the reader's
convenience we shall include a brief proof, following Nishikawa\ \cite{N}%
.\medskip

\noindent \textbf{Lemma. }\textit{Let }$M^{n}$\textit{\ be an }$n$\textit{%
-dimensional, }$n\geq 2$,\textit{\ complete Riemannian manifold with Ricci
curvature }$Ric\geq -\left( n-1\right) k^{2}$,\textit{\ where }$k$ \textit{%
is a positive constant}.\textit{\ Suppose that }$g$\textit{\ is a smooth
non-negative function on }$M^{n}$ \textit{satisfying} 
\begin{equation}
\Delta g\geq cg^{2},  \label{3.1}
\end{equation}
\textit{where }$c$ \textit{is\ a positive constant and }$\Delta $ \textit{%
stands for the Laplacian operator. Then }$g$ \textit{vanishes identically.}%
\medskip

\noindent \textit{Proof.} Take a point $x_{0}\in M^{n}$. If $g$ attains its
maximum at $x_{0}$ then $g\left( x_{0}\right) =0$ and the result is true. We
shall prove that $g\left( x\right) \leq g\left( x_{0}\right) ,$ for any
point $x\in M^{3}$. Suppose in the contrary, that\ there exists a point $%
x_{1}$ such that $g\left( x_{0}\right) <g\left( x_{1}\right) $. We\ set $%
\beta =d\left( x_{0},x_{1}\right) $ and denote by $B_{a}\left( x_{0}\right) $
the geodesic ball of radius $a$ centered at $x_{0}$, where $d$ is the
distance on $M^{3}$. We consider the function $G:B_{a}\left( x_{0}\right)
\rightarrow \mathbb{R}$, $G\left( x\right) =\left( a^{2}-\rho ^{2}\left(
x\right) \right) ^{2}g\left( x\right) $, where $a>\beta $ and $\rho \left(
x\right) =d\left( x,x_{0}\right) $. Since $G\left| _{\partial B_{a}\left(
x_{0}\right) }\right. =0$ and $G$ is non-negative we deduce that $G$ attains
its maximum at some point $x_{2}\in B_{a}\left( x_{0}\right) $. First let us
assume that $x_{2}$ is not a cut point of $x_{0}$. Then $G$ is smooth near $%
x_{2}$ and by the maximum principle we have 
\begin{equation}
\nabla G\left( x_{2}\right) =0,\quad \Delta G\left( x_{2}\right) \leq 0.
\end{equation}
It is well known that $\rho \Delta \rho \left( x_{2}\right) \leq \left(
n-1\right) \left( 1+k\rho \left( x_{2}\right) \right) $ since the Ricci
curvature is bounded from below by $-\left( n-1\right) k^{2}$. Then from $(4)$, we get 
\begin{equation}
G\left( x_{2}\right) \leq \frac{a^{2}}{c}\left( 28+4k\left( n-1\right)
a\right) .  \label{3.3}
\end{equation}
Assume now that $x_{2}$ is a cut point of $x_{0}$. Let $\sigma $ be a unit
speed segment from $x_{0}$ to $x_{2}$ and $\overline{x}_{0}\in \sigma $ such
that $d\left( x_{0},\overline{x}_{0}\right) =\varepsilon $, for sufficiently
small $\varepsilon >0$. Then the function $\overline{\rho }\left( x\right)
=d\left( x,\overline{x}_{0}\right) $ is smooth near $x_{2}$ and the function 
\begin{equation*}
G_{\varepsilon }\left( x\right) =\left( a^{2}-\left( \overline{\rho }\left(
x\right) +\varepsilon \right) ^{2}\right) g\left( x\right)
\end{equation*}
is a ``support function'' of $G$, i.e. $G_{\varepsilon }\left( x\right) \leq
G\left( x\right) $, $G_{\varepsilon }\left( x_{2}\right) =G\left(
x_{2}\right) $. Thus $G_{\varepsilon }$ attains its maximum at $x_{2}$.
Proceeding as before and passing to the limit we get the same estimate as in 
$(5)$. Since $x_{2}$ is the maximum point of $G$ in $%
B_{a}\left( x_{0}\right) $ we have $G\left( x_{1}\right) \leq G\left(
x_{2}\right) $ and thus 
\begin{equation}
g\left( x_{1}\right) \leq \frac{a^{2}\left( 28+4k\left( n-1\right) a\right) 
}{c\left( a^{2}-\rho ^{2}\left( x_{1}\right) \right) ^{2}}.  \label{3.4}
\end{equation}
Since $M^{n}$ is complete, letting $a\rightarrow \infty $, from $(6)$, we get $g\left( x_{1}\right) =0$ which is a
contradiction.\medskip

\textit{Proof of Theorem }1. Without loss of generality we may assume that $%
M^{3}$ is simply connected, after passing to the universal covering space.
Since $M^{3}$ is simply connected, the standard monodromy argument allows us
to define a global orthonormal frame field $\left\{
e_{1},e_{2},e_{3}\right\} $ of principal directions. The assumptions of
Theorem 1 imply that $M^{3}$ has three distinct principal curvatures $%
\lambda >0>-\lambda $, and Ricci curvature bounded from below. The functions 
$u$ and $v$ are well defined on entire $M^{3}$. Moreover, we claim that $u$
and $v$ are harmonic functions. Indeed, 
\begin{eqnarray*}
\Delta v &=&e_{1}e_{1}\left( v\right) +e_{2}e_{2}\left( v\right)
+e_{3}e_{3}\left( v\right) \\
&&-\left( \omega _{31}\left( e_{3}\right) +\omega _{21}\left( e_{2}\right)
\right) e_{1}\left( v\right) \\
&&-\left( \omega _{12}\left( e_{1}\right) +\omega _{32}\left( e_{3}\right)
\right) e_{2}\left( v\right) \\
&&-\left( \omega _{13}\left( e_{1}\right) +\omega _{23}\left( e_{2}\right)
\right) e_{3}\left( v\right) \\
&=&e_{1}e_{1}\left( v\right) +e_{2}e_{2}\left( v\right) +e_{3}e_{3}\left(
v\right) \\
&&-\frac{1}{2}e_{1}\left( \log \lambda \right) e_{1}\left( v\right)
-2ve_{2}\left( v\right) \\
&&-\frac{1}{2}e_{3}\left( \log \lambda \right) e_{3}\left( v\right) .
\end{eqnarray*}
Making use of $(1)$, we get 
\begin{eqnarray*}
e_{1}e_{1}\left( v\right) &=&-e_{1}e_{3}\left( u\right) , \\
e_{2}e_{2}\left( v\right) &=&2ve_{2}\left( v\right) -2ue_{2}\left( u\right)
\\
&=&2v^{3}-6vu^{2}, \\
e_{3}e_{3}\left( v\right) &=&e_{3}e_{1}\left( u\right) .
\end{eqnarray*}
Therefore, taking $(2)$ into account, we obtain 
\begin{eqnarray*}
\Delta v &=&-e_{1}e_{3}\left( u\right) +e_{3}e_{1}\left( u\right)
+2v^{3}-6vu^{2} \\
&&-\frac{1}{2}e_{1}\left( \log \lambda \right) e_{1}\left( v\right)
-2ve_{2}\left( v\right) -\frac{1}{2}e_{3}\left( \log \lambda \right)
e_{3}\left( v\right) \\
&=&\frac{1}{2}e_{3}\left( \log \lambda \right) e_{1}\left( u\right)
+2ue_{2}\left( u\right) -\frac{1}{2}e_{1}\left( \log \lambda \right)
e_{3}\left( u\right) \\
&&+2v^{3}-6vu^{2} \\
&&-\frac{1}{2}e_{1}\left( \log \lambda \right) e_{1}\left( v\right)
-2ve_{2}\left( v\right) -\frac{1}{2}e_{3}\left( \log \lambda \right)
e_{3}\left( v\right) \\
&=&0.
\end{eqnarray*}
In a similar way, we verify that $\Delta u=0$.

Using $(1) $, and the fact that $\Delta v=\Delta u=0$ we get 
\begin{eqnarray*}
\frac{1}{2}\Delta \left( u^{2}+v^{2}\right) &=&\left\| \nabla u\right\|
^{2}+\left\| \nabla v\right\| ^{2} \\
&\geq &\left( e_{2}\left( u\right) \right) ^{2}+\left( e_{2}\left( v\right)
\right) ^{2} \\
&=&4u^{2}v^{2}+\left( v^{2}-u^{2}\right) ^{2}.
\end{eqnarray*}
Thus, 
\begin{equation*}
\Delta \left( u^{2}+v^{2}\right) \geq 2\left( u^{2}+v^{2}\right) ^{2}.
\end{equation*}
Appealing to the Lemma, we infer that $u^{2}+v^{2}$ is identically zero.
Thus $u=v\equiv 0$ and $\lambda $ is constant along the integral curves of $%
e_{2}$. Consider the $2$-dimensional distribution $V$ which is spanned by $%
e_{1}$ and $e_{3}$. Because $u\equiv 0$, from $(2)$ we see
that $V$ is involutive. Let $L_{x}^{2}$ be\ a maximal integral submanifold
of $V$ passing through a point $x$ of $M^{3}$, and denote by $%
i:L_{x}^{2}\rightarrow M^{3}$ its inclusion map. Then $\widetilde{f}=f\circ
i:L_{x}^{2}\rightarrow \mathbb{R}^{4}$ defines an immersion. Let $A_{1}$ and 
$A_{2}$ be the shape operators of $\widetilde{f}$ in the directions $%
df\left( e_{2}\right) $ and $\xi $, respectively. A direct calculation shows
that $A_{1}=0$ and 
\begin{equation*}
A_{2}\sim \left( 
\begin{array}{cc}
\lambda  & 0 \\ 
0 & -\lambda 
\end{array}
\right) 
\end{equation*}
with respect to the basis $\left\{ e_{1},e_{3}\right\} $. Thus, $\widetilde{f%
}:L_{x}^{2}\rightarrow \mathbb{R}^{4}$ is a minimal surface with bounded
Gaussian curvature which lies in $\mathbb{R}^{3}$ and $df\left( e_{2}\right) 
$ is constant along $\widetilde{f}$. Hence $f\left( M^{3}\right) $ splits as
we wished, and\ this completes the proof of theorem.\medskip

\textit{Proof of Theorem}\textbf{\ }2\textbf{. }We can choose an orientation
such that $k_{1}>0>k_{2}\geq k_{3}$. According to the Principal Curvature
Theorem due to Smyth and Xavier \cite{SX}, we have $\inf \left( \Lambda \cap 
\mathbb{R}^{+}\right) =0$, where $\Lambda $ is the set of values assumed by\
the non-zero principal curvatures. Hence there exists a secuence of points $%
\left\{ x_{n}\right\} $ such that $k_{1}\left( x_{n}\right) \rightarrow 0$.
From the minimality we get $k_{1}\geq \left| k_{i}\right| $, $i=2,3$ and
consequently the Gauss-Kronecker curvature satisfies $K\left( x_{n}\right)
\rightarrow 0$. Thus $K$ is zero, and this completes the proof.

\end{document}